\documentclass[12pt]{article}
\usepackage{amsfonts}
\usepackage{}
\usepackage{wasysym,amssymb,eufrak,indentfirst,graphicx,cite,amsthm,color}
\usepackage[bookmarksnumbered, colorlinks, plainpages]{hyperref}
\newtheorem{Theorem}{Theorem}[section]
\theoremstyle{definition}

\newtheorem{Lemma}[Theorem]{Lemma}

\theoremstyle{definition}

\def\nd{\mathrel{\bigm|\kern-.7em/}}
\def\f{\noindent}

\def\demo{\f{ \textbf{Proof}}\hskip10pt}
\def\qed{\hfill $\Box$}

\begin{document}
\baselineskip 17pt

\title{The permutability of $\sigma_i$-sylowizers of some $\sigma_i$-subgroups in finite
groups}

\author{Zhenya Liu AND Wenbin Guo}
\date{}
\maketitle

\begin{abstract}
 Let
 $\sigma=\{\sigma_{i}|i\in I\}$ be a partition of the set of all primes $\mathbb{P}$, $G$ a finite group and
 $\sigma(G)=\{\sigma_{i}|\sigma_{i}\cap \pi(|G|)\neq\emptyset\}$.
  A subgroup $S$ of a group $G$ is called a $\sigma_i$-sylowizer of a
  $\sigma_i$-subgroup $R$ in $G$ if $S$ is maximal in $G$ with respect to having $R$ as its
Hall $\sigma_i$-subgroup. The main aim of this paper is to investigate the
influence of $\sigma_i$-sylowizers on the structure of finite groups. We obtained some
new characterizations of supersoluble groups by the
permutability of the $\sigma_i$-sylowizers of some $\sigma_i$-subgroups.
\end{abstract}

\let\thefootnoteorig\thefootnote
\renewcommand{\thefootnote}{\empty}


 \let\thefootnote\thefootnoteorig

\section{Introduction}

Let $\pi$ denotes a set of primes. The concept  of $\pi$-Sylowizers has been introduced by W. Gaschutz\cite{WG}.
If $R$ is a $\pi$-subgroup of the group $G$, then a $\pi$-Sylowizer of $R$ in $G$ is a subgroup $S$ of $G$ maximal with respect to containing $R$ as a Hall $\pi$-subgroup.

$\mathbb{P}$ is the set of all primes and $n$ is a natural number. Let $\sigma=\{\sigma_i |i\in I\}$ is some partition of all primes $\mathbb{P}$, that is, $\mathbb{P}=\bigcup_{i\in I}\sigma_i$ and $\sigma_i\cap \sigma_j=\emptyset$ for all $i\neq j$.
We write $\sigma(G)=\{\sigma_i|\sigma_i\cap \pi(G)\neq \emptyset\}$.

Following \cite{GSS}, two subgroups $H$ and $T$ of
a group $G$ are conditionally permutable (or in brevity, $c$-permutable) in $G$ if
there exists an element $x \in G$ such that $HT^x = T^x H$.

\section{Preliminaries}

\begin{Lemma}
Let $H$ be a $\sigma_i$-subgroup of $G$ for some $\sigma_i\in \sigma(G)$. Assume that $K$ is a subgroup satisfying $H \leq  K \leq G$ and $T$ is a $\sigma_i$-sylowizer
 of $H$ in $K$. Then there is a $\sigma_i$-sylowizer
$S$ of $H$ in $G$ such that $T = S \cap K$.
\end{Lemma}

\demo
 Since $H$ is a Hall $\sigma_i$-subgroup of $T$, there is a $\sigma_i$-sylowizer $S$ of $H$ in $G$
such that $S \geq T$. Then $H$ is a Hall $\sigma_i$-subgroup of $S \cap K$. Since $T \leq  S \cap K$
and $T$ is a $\sigma_i$-sylowizer of $H$ in $K$, we get $T = S \cap K$ by the maximality
of $T$.
\qed

\begin{Lemma}
Let $R$ be a $\sigma_i$-subgroup of $G$ for some $\sigma_i\in \sigma(G)$. Assume that $N$ is a normal subgroup
of $G$ and $R$ is a Hall $\sigma_i$-subgroup of $RN$. Then $S$ is a $\sigma_i$-sylowizer of $R$ in $G$
if and only if $S/N$ is a $\sigma_i$-sylowizer of $RN/N$ in $G/N$.
\end{Lemma}

\demo
 Let $S$ be a $\sigma_i$-sylowizer of $R$ in $G$. Since $R$ is a Hall $\sigma_i$-subgroup of
$RN$, $R$ is a Hall $\sigma_i$-subgroup of $SN$. Thus $N\leq S$ by the maximality of
$S$ and so $RN/N$ is a Hall $\sigma_i$-subgroup of $S/N$. If $S/N$ is not a $\sigma_i$-sylowizer
of $RN/N$ in $G/N$, then there is a $\sigma_i$-sylowizer $S_0 /N$ of $RN/N$ in $G/N$ such
that $S_0 /N > S/N$. Now, $S_0 > S$ and $R$ is a Hall $\sigma_i$-subgroup of $S_0$, which
contradicts the fact that $S$ is a $\sigma_i$-sylowizer of $R$ in $G$. Thus $S/N$ is a $\sigma_i$-sylowizer
of $RN/N$ in $G/N$.

Conversely, if $S/N$ is a $\sigma_i$-sylowizer of $RN/N$ in $G/N$, then $R$ is a Hall
$\sigma_i$-subgroup of $S$. If $S$ is not a $\sigma_i$-sylowizer of $R$ in $G$, then there is a $\sigma_i$-sylowizer
$S_0$ of $R$ in $G$ such that $S_0 > S$. Therefore $RN/N$ is a Hall $\sigma_i$-subgroup of
$S_0 /N$, which contradicts the fact that $S/N$ is a $\sigma_i$-sylowizer of $RN/N$ in $G/N$.
Thus $S$ is a $\sigma_i$-sylowizer of $R$ in $G$.
\qed

\begin{Lemma}
Let $R$ be a $\sigma_i$-subgroup of a $\sigma$-full group $G$ for some $\sigma_i\in \sigma(G)$ and $S$ a $\sigma_i$-sylowizer of $R$ in $G$. If $S$ is
$\sigma$-permutable in $G$, then $O^{\sigma_i} (G) \leq S$. In particular, $S = RO^{\sigma_i}(G)$ is the unique
$\sigma_i$-sylowizer of $R$ in $G$.
\end{Lemma}

\demo
Let Q be a Hall $\sigma_j$-subgroup of $G$ with $\sigma_j\in \sigma(G)$ and $\sigma_i\cap \sigma_j=\emptyset$. Since $S$ is
$\sigma$-permutable, we have $SQ \leq G$. Note that since $R$ is a Hall $\sigma_i$-subgroup of
$SQ$, we have $QS = S$ by the maximality of $S$. Hence $Q \leq S$. It shows that
$O^{\sigma_i} (G) \leq S$.
\qed

\begin{Lemma}
Let $R$ be a $\sigma_i$-subgroup of a $\sigma$-full group of Sylow type $G$ for some $\sigma_i\in \sigma(G)$ and $S$ a $\sigma_i$-sylowizer of $R$ in $G$. Then
$S$ is $c$-permutable with every Hall $\sigma_j$-subgroup of $G$ for all $\sigma_j\in \sigma(G)$ if and only if $|G : S|$ is a $\sigma_i$-number.
\end{Lemma}

\demo
The sufficiency is evident, we only need to prove the necessity.

Let $Q$ be a Hall $\sigma_j$-subgroup of $G$ with $\sigma_j\in \sigma(G)$ and $\sigma_i\cap \sigma_j=\emptyset$. Since $S$ is $c$-permutable
with $Q$, we have $SQ^x = Q^x S$ for some element $x \in G$. Note that since $R$ is a
Hall $\sigma_i$-subgroup of $SQ^x$, we have $Q^x S = S$ by the maximality of $S$. Hence
$Q^x \leq S$. It implies that $|G : S|$ is a $\sigma_i$-number.
\qed

\begin{Theorem}
Let $G$ be a $\sigma$-full group of Sylow type and $\mathcal{H}=\{H_1,\cdots,H_t\}$ be a complete Hall $\sigma$-set of $G$ such that $H_i$ is a
nilpotent $\sigma_i$-subgroup for all $i=1,\cdots,t$. Suppose that for any $\sigma_i\in \sigma(G)$, every
maximal subgroup of any non-cyclic $H_i$ has a $\sigma_i$-sylowizer that is
$c$-permutable with every member of $\mathcal{H}$, then $G$ is supersoluble.
\end{Theorem}

\demo
Assume that this is false and let $G$ be a counterexample of minimal order.
Then:

$(1)$ {\sl Let $N$ be a minimal normal subgroup of $G$. Then $G$ is supersoluble}.

We consider the quotient group $G/N$. It is clear that $G/N$ is a $\sigma$-full group of Sylow type and $\mathcal{H}N/N$ is a complete Hall $\sigma$-set of $G/N$ such that $H_iN/N$ is nilpotent.
Let $H/N$ be a maximal subgroup of $H_iN/N$ and $H_{\sigma_i}$ be a Hall $\sigma_i$-subgroup of $H$ contained in $H_i$. Then $H=H_{\sigma_i}N$.
Since $H_{\sigma_i}\cap N=N_{\sigma_i}=H_i\cap N$, where $N_{\sigma_i}$ denotes a Hall $\sigma_i$-subgroup of $N$, we have that
$$|H_i:H_{\sigma_i}|=\frac{|H_i||N|}{|H_i\cap N|}\cdot \frac{|H_{\sigma_i}\cap N|}{|H_{\sigma_i}||N|}=|H_iN:H|=q$$
for some $q\in \sigma_i$. This shows that $H_{\sigma_i}$ is a maximal subgroup of $H_i$.
If $H_iN/N$ is non-cyclic, then so is $H_i$. Thus if $S/N$ is a $\sigma_i$-sylowizer
of $H/N$ in $G/N$, then $S$ is a $\sigma_i$-sylowizer of $H_{\sigma_i}$ in $G$ by Lemma 2.2. Moreover, if
$S$ is $c$-permutable with every member of $\mathcal{H}$, then $S/N$ is $c$-permutable
with every member of $\mathcal{H}N/N$ by Lemma 2.4. It shows that $G/N$ satisfies
the hypotheses. Thus $G/N$ is supersoluble by the choice of $G$.

$(2)$ {\sl $N$ is the unique proper minimal normal subgroup of $G$ and $\Phi(G)=1$}.

Let $p$ be the smallest prime divisor of $G$ and $p\in \sigma_i$.
If $H_i$ is cyclic, then $G$ is $p$-nilpotent.
This shows that $G$ has a
proper minimal normal subgroup.
Thus we may assume that $H_i$
is non-cyclic. Let $M$ be a maximal subgroup of $H_i$ of index $p$ and $S$
a $\sigma_i$-sylowizer of $M$ in $G$ that is
$c$-permutable with every member of $\mathcal{H}$.
Then $|G : S| = p$ by Lemma 2.4 and so $S \unlhd G$. Therefore we may choose a
proper minimal normal subgroup of $G$ contained in $S$, say $N$. By Claim $(1)$, $G/N$
is supersoluble. Moreover, $N$ is the unique minimal normal subgroup of $G$.
Since the class of all supersoluble groups is a saturated formation, we may
assume further that $|\Phi(G) |= 1$.

$(3)$ {\sl $N$ is soluble}.

Assume that $N$ is not soluble. Then $p=2$ and $2 | |N|$.
Let $P$ be a Sylow $2$-subgroup of $H_i$. Then $N_2=P\cap N$ is a Sylow $2$-subgroup of $N$. If $N_2 \leq \Phi(H_i)$,
then $N_2\leq \Phi(P)$, and so $N$ is $2$-nilpotent by Tate's theorem, a contradiction. Hence $N_2 \nleq  \Phi(H_i)$.
Thus there is a maximal subgroup $K$ of $H_i$ such that $H_i = KN_2 $. Let $S_0$ be a
$\sigma_i$-sylowizer of $K$ in $G$ that is $c$-permutable with every member of $\mathcal{H}$.
Then $|G : S_0 | = 2$ by Lemma 2.4. Thus $G = S_0 H_i = S_0 N_2 = S_0 N$. Now,
$|N : N \cap S_0 | = |G : S_0 | = 2$,
which implies that $N \cap S_0 \unlhd N$. Since $N \cap S_0 \unlhd S_0$, we have $N \cap S_0 \unlhd G$. Note
that $N$ is a minimal normal subgroup of $G$, we have $N \cap S_0 = 1$. Thus
$|N| = |G : S_0 | = 2$, a contradiction.

$(4)$ {\sl Final contradiction}.

By Claim $(3)$, we may assume that $N$ is a $q$-subgroup for some prime $q\in \sigma_j$. Since
$\Phi(G) = 1$, there is a maximal subgroup $T$ of $G$ such that $G = TN$. Let $T_{\sigma_j}$
be a Hall $\sigma_j$-subgroup of $T$ contained in $H_j$. Then $H_j = T_{\sigma_j} N$ is a Hall $\sigma_j$-subgroup of $G$.
If $H_j$ is cyclic, then $G$ is supersoluble by the supersolublity of $G/N$. Thus
we may assume that $H_j$ is non-cyclic. Let $Q \geq T_{\sigma_j}$ be a maximal subgroup
of $H_j$ and $Y$ a $\sigma_j$-sylowizer of $Q$ in $G$ that is $c$-permutable with every member
of $\mathcal{H}$. Then $|G : Y | = q$ by Lemma 2.4 and $N \nleq Y  $. Otherwise
$H_j = Q N \leq  Y $, which contradicts the fact that $Q$ is a Hall $\sigma_j$-subgroup of $Y$. Thus $G = Y N$
and so $|N| = |G : Y | = q$. It implies that $G$ is supersoluble, a contradiction.
This contradiction completes the proof.
\qed

\begin{Theorem}
Let $\mathfrak{F}$ be a soluble saturated formation containing all supersoluble groups
and let $E$ be a normal subgroup
 of $G$ with $G/E\in \mathfrak{F}$. Suppose that $G$ is a $\sigma$-full group of Sylow type
 and $\mathcal{H}=\{H_1,\cdots,H_t\}$ is a complete Hall $\sigma$-set of $G$ such that $H_i$ is a
nilpotent $\sigma_i$-subgroup for all $i=1,\cdots,t$.
If for any $\sigma_i\in \sigma(E)$, every maximal subgroup of any non-cyclic $H_i\cap E$ has a $\sigma_i$-sylowizer that is
$c$-permutable with every member of $\mathcal{H}$, then $G\in \mathfrak{F}$.
\end{Theorem}

\demo
The conclusion holds when $E = G$ by Theorem 2.5, thus we may assume
that $E < G$. Let $N$ be a minimal normal subgroup of $G$ contained in $E$.

$(1)$ {\sl $E$ is supersoluble.}

Let $Q$ be a maximal subgroup of a non-cyclic Hall $\sigma_i$-subgroup $H_i\cap E$ of $E$ and
$S$ a $\sigma_i$-sylowizer of $Q$ in $G$ that is $c$-permutable with member of $\mathcal{H}$.
By Lemma 2.4, $|G : S |$ is a $\sigma_i$-number.
Let $Y = S  \cap E$. Since
$|E : Y  | = |E : S \cap E| = |S  E : S  |$ divides $|G : S  |$, $|E : Y  |$ is a $\sigma_i$-number.
Hence $Y $ is a $\sigma_i$-sylowizer of $Q$ in $E$ and $Y $ is $c$-permutable with every member of
$\mathcal{H}\cap E$ by Lemma 2.4. Thus $E$ is supersoluble by Theorem 2.5.

$(2)$ {\sl $N$ is the unique minimal normal subgroup of $G$ contained in $E$ and
$N \cap \Phi(G) = 1$.}

Consider the quotient group $G/N$, evidently $(G/N)/(E/N) \in  \mathfrak{F}$.
Since $E$ is supersoluble by Claim $(1)$, we have that $N$ is a $p$-group for some prime $p$. Without loss of generality, we may write
$E_i=H_i \cap E$ for all $i\in \{1,\cdots,t\}$ and assume that $p\in \sigma_i$ for some $i$. Let $J/N$ be a maximal
subgroup of $E_i/N$, then $J$ is a maximal subgroup of $E_i$.
If $S/N$ is a $\sigma_i$-sylowizer of $J/N$ in $G/N$, then $S$ is a $\sigma_i$-sylowizer
of $J$ in $G$ by Lemma 2.2. Moreover, if $S$ is $c$-permutable with every member of $\mathcal{H}$,
then $S/N$ is $c$-permutable with every member of $\mathcal{H}N/N$
by Lemma 2.4.
Let $J/N$ be a maximal
subgroup of $E_jN/N$ and $J_{\sigma_j}$ a Hall ${\sigma_j}$-subgroup of $J$ contained in $E_j$, where $i\neq j$. Then $J_{\sigma_j}$ is a maximal subgroup of $E_j$.
If $S/N$ is a $\sigma_j$-sylowizer of $J_{\sigma_j}N/N$ in $G/N$, then $S$ is a $\sigma_j$-sylowizer
of $J_{\sigma_j}$ in $G$ by Lemma 2.2. Moreover, if $S$ is $c$-permutable with every member of $\mathcal{H}$,
then $S/N$ is $c$-permutable with every member of $\mathcal{H}N/N$
by Lemma 2.4.
This shows that $(G/N,E/N)$ satisfies the hypotheses. Thus we
may have that $G/N \in  \mathfrak{F}$ by induction. Moreover, $N$ is the unique minimal
normal subgroup of $G$ contained in $E$ and $N \cap  \Phi(G) = 1$.

$(3)$ {\sl $N$ is an elementary abelian $p$-subgroup, where $p$ is the largest prime
divisor of $|E|$.}

Since $E$ is supersoluble by Claim $(1)$, the Sylow $p$-subgroup $E_P$ of $E$ is normal in $G$. Note that $N$ is the unique minimal
normal subgroup of $G$ contained in $E$, $N \leq  E_P$ is an elementary abelian
$p$-subgroup.

$(4)$ {\sl $G \in  \mathfrak{F}$.}

Without loss of generality, we may assume that $p\in \sigma_i$.
If $E_i$ is cyclic, then $|N| = p$ and so $G \in  \mathfrak{F}$.
 Assume that
$E_i$ is non-cyclic. Since $N\nleq \Phi(G)$, there is a maximal subgroup $M$ of $G$ such
that $G = MN$ and $M \cap N = 1$. Thus $E_i = N(M \cap E_i )$ and $H_i= NM \cap H_i = N(M \cap H_i) = NM_i$.
Since $M_i < H_i $, we may choose $P \lessdot H_i$ such that $M_i \leq P $. Since $M\cap E_i \leq P $,
$P\cap E_i = P \cap N(M \cap E_i ) = (P \cap N)(M \cap E_i)$. Note that $M \cap N = 1$,
we have
$$|E_i : E_i \cap P| = |N(M \cap E_i) : (P\cap  N)(M \cap E_i )| = |N : P \cap N| = p.$$
Hence $R = E_i \cap P$ is a maximal subgroup of $E_i$. Let $S$ be a $\sigma_i$-sylowizer of $R$
in $G$ that is $c$-permutable with every member of $\mathcal{H}$. Then $|G : S|$ is a $\sigma_i$-number
by Lemma 2.4. Since $G$ is soluble, we may write $S = RS_{\sigma_i'}$ and $M = M_iM_{\sigma_i'}$.
Note also that $|G:S|$ and $|G:M|$ are $\sigma_i$-number in $G$, $S_{\sigma_i'}$ and $M_{\sigma_i'}$ are also
Hall ${\sigma_i'}$-subgroups of $G$. Thus there is an element $g$ of $G$ such that
$S_{\sigma_i'}^g= M_{\sigma_i'}$.
Since $G = H_iS^g $, we may write $g = xy$, where $x \in H_i$ and $y \in S^g$. Note that
since $R = E_i \cap P \unlhd H_i$, we have $R^y = R^{xy} \leq S^g$ and so $R \leq S^g$. Thus
$S^g = RM_{\sigma_i'}$. Since $RM_i= (P \cap E_i )M_i = P \cap E_i M_i = P \cap NM_i = P \leq G$,
we have $RM \leq G$. Since $M$ is a maximal subgroup, either $RM = M$ or
$RM = G$.

If $RM = G$, then $RM_i= P$ is a Hall $\sigma_i$-subgroup of $G$, which is impossible.
Thus $RM = M$ and so $R \leq M\cap E_i$. Since $G = MN = ME_i$, we have $E_i\nleq M$.
Note that since $R \lessdot E_i$, we have $R = M \cap E_i$. Thus
$|N| = |G : M| = |E_i :E_i \cap M| = |E_i : R| = p$. By \cite[Theorem 2]{BP}, $G \in  \mathfrak{F}$, as required.
\qed

\end{document}